\documentclass[12pt]{amsart}
\usepackage[letterpaper,hmargin=1.1in,vmargin=1.05in]{geometry}
\usepackage{amsmath,amssymb,amsthm,mathtools,amscd}
\usepackage{tikz-cd}
\usepackage[T1]{fontenc}
\usepackage{lmodern,microtype,needspace}
\usepackage[hidelinks]{hyperref}
\newcommand{\PP}{\mathbf P}
\renewcommand{\AA}{\mathbf A}
\newcommand{\CC}{\mathbf C}
\newcommand{\OO}{\mathcal O}
\newcommand{\sm}{\mathrm{sm}}
\DeclareMathOperator{\Bl}{Bl}
\DeclareMathOperator{\Sing}{Sing}
\DeclareMathOperator{\pr}{pr}
\newtheorem{theorem}{Theorem}[section]
\newtheorem{proposition}[theorem]{Proposition}
\newtheorem{lemma}[theorem]{Lemma}
\newtheorem{corollary}[theorem]{Corollary}
\title{{Algebraic ellipticity of double covers of
$(\mathbf P^1)^3$}}
\author{Alexander Dvorsky}
\address{Department of Mathematics, University of Miami, 
Coral Gables, FL 33124, USA}
\email{dvorsky@gmail.com}
\author{Shulim Kaliman}
\address{Department of Mathematics, University of Miami, 
Coral Gables, FL 33124, USA}
\email{kaliman@math.miami.edu}
\author{Mikhail Zaidenberg}
\address{Univ. Grenoble Alpes, CNRS, IF, 
38000 Grenoble, France}
\email{mikhail.zaidenberg@gmail.com}
\date{}
\expandafter\providecommand\csname subjclassname@2020\endcsname{2020 Mathematics Subject Classification}
\subjclass[2020]{Primary 14J45; Secondary 14D06, 14M20, 32Q56}
\keywords{Algebraic Gromov ellipticity; Fano threefold; del Pezzo fibration; algebraic spray; conic bundle}
\begin{document}

\begin{abstract}
{We show that every smooth complex double cover of $(\mathbf P^1)^3$ branched along a smooth divisor of tridegree $(2,2,2)$ is algebraically Gromov elliptic. These irrational Fano threefolds form Mori--Mukai family 3-1 and have anticanonical genus seven and Picard rank three. Their three coordinate projections define degree-four del Pezzo fibrations. Each threefold also admits a dominant rational map of degree two from $\mathbf P^3$.}
\end{abstract}
\maketitle

\section{Introduction}
All varieties considered in this paper are defined over $\CC$. 
Recall that an algebraic spray on a smooth variety $V$ 
is a vector bundle $p:E\to V$ with a morphism $s:E\to V$ 
satisfying $s(0_v)=v$. It is dominating at $v$ if the differential 
of $s|_{E_v}$ at $0_v$ surjects onto $T_vV$. The variety $V$ 
is said to be algebraically Gromov elliptic, or simply  algebraically elliptic, 
if it admits a spray dominating at every point. 
According to Banecki \cite{Ban},
every smooth complete rational variety is algebraically elliptic. 
The first examples of irrational smooth
projective {algebraically elliptic} varieties are smooth cubic threefolds and
{very general smooth cubic fourfolds, see \cite[Theorem 1.1]{KZCub}} 
and the introduction in \cite{DKZ}.
Very general smooth quartic double spaces of dimension 3 and 4 are also 
{algebraically elliptic and are not stably rational}, see  \cite{DKZ}. 
In this paper we extend this list of examples by the following irrational Fano threefolds.

{For a smooth effective divisor $R$ on $Y=(\PP^1)^3$ of tridegree $(2,2,2)$, let $\pi:X\to Y$ be the double cover branched along $R$. Its three coordinate fibrations are $f_i=\pr_i\circ\pi:X\to\PP^1$, $i=1,2,3$.}

\begin{theorem}\label{thm:main}
For every smooth $R\in|\OO_Y(2,2,2)|$, the threefold $X$ is algebraically Gromov elliptic.
\end{theorem}

These threefolds form family 3-1 in the Mori--Mukai classification \cite{MM} and have Picard number three. Let $L=\OO_Y(1,1,1)$. Since $K_Y=-2L$ and $L^3=6$, the double-cover formula $K_X=\pi^*(K_Y+L)$ gives
\begin{equation}\label{eq:fano}
 -K_X=\pi^*(L),\qquad (-K_X)^3=2\cdot6=12.
\end{equation}
{Thus $X$ is Fano of anticanonical genus} $1+\frac12(-K_X)^3=7$. According to  Prokhorov  \cite[Theorem 6.1, case (6.1.g)]{Prokhorov}, {every smooth threefold} $X$ as above is irrational, whereas it is unirational. Moreover, we provide in Proposition~\ref{prop:unirational} a dominant rational map $\PP^3\dashrightarrow X$ of degree two.

{The proof uses local rank-one sprays. For a smooth variety $V$, such a spray is a morphism}
\[
 s:U\times\AA^1\longrightarrow V,\qquad s(u,0)=u,
\]
where $U\subset V$ is Zariski open; its values need not remain in $U$. {Its velocity at $a\in U$ is $\partial_t s(a,0)\in T_aV$.} To {establish} the algebraic ellipticity of $V$, it suffices to construct, near every point of $V$, finitely many {local rank-one sprays} whose velocities span its tangent space. {By Gromov's localization lemma, these local sprays extend to global rank-one sprays whose velocity directions at the chosen point are unchanged}, see  \cite[Lemma 1.3 and the proof of Lemma 1.5]{KZ}. 
{Openness of domination and quasicompactness give finitely many global sprays whose velocities span the tangent space at every point of $V$. Thus $V$ is algebraically subelliptic, and the equivalence of subellipticity and ellipticity \cite[Theorem 0.1]{KZ} gives algebraic ellipticity.}  

We obtain these local sprays on our threefold $X$ from the three natural del Pezzo fibrations $f_i:X\to\PP^1$. Their fibers are degree-four del Pezzo surfaces with at worst Du Val singularities, and their general fibers are smooth (Lemma~\ref{lem:fibers}). {At a smooth point $a$ of a fiber $S$, which may be singular elsewhere, Lemma~\ref{lem:centers} provides} an auxiliary point $q\ne a$, lying on no line of $S$, such that the anticanonical section through $a$ and singular at $q$ is integral with a single ordinary node. Blowing up $q$ gives a cubic del Pezzo surface. {If $H$ denotes the pullback of $-K_S$ to this blow-up and $e_q$ its exceptional curve, the anticanonical class is $H-e_q$, whereas} the strict transforms of the sections singular at $q$ form a conic pencil of class $H-2e_q$. The chosen nodal member becomes a smooth conic.

To carry out this construction in nearby surface fibers, we extend $q$ to a section of $f_i$ and blow up that section. This extension uses weak approximation, including when $S$ is singular. 
{We pull the smooth conic family back by the map that assigns to each nearby point $u$ the member of the pencil through $u$. The point $u$ itself gives a section of the pullback, which is therefore a $\PP^1$-bundle, trivial on a suitable Zariski neighborhood $U$ of $a$.} An affine chart with the marked section at zero, followed by the blow-down to $X$, gives a local rank-one spray on $U$ with values in $X$.

{We choose the auxiliary point twice so that the corresponding curves have different tangent directions at $a$. The two spray velocities then span $T_aS$.} Finally, at every point of $X$, at least two of the three fibrations $f_i$ are smooth there and have distinct fiber tangent planes. The corresponding local sprays span $T_aX$. 
{This proves Theorem~\ref{thm:main}.}

\paragraph{{AI disclosure.}}
{Motivated by their work on algebraic ellipticity of quartic double spaces \cite{DKZ}, the authors directed GPT-6 Astra to search for further examples of nonrational algebraically elliptic varieties. This search led to the class considered in this paper. The resulting arguments were revised, expanded where necessary, and streamlined by the authors with assistance from GPT-6 Astra. The authors take responsibility for the mathematical content.}

\section{Anticanonical curves through a prescribed point}
A weak del Pezzo surface is a smooth projective surface whose anticanonical divisor is nef and big. In degrees three and four, its anticanonical system is base point free and defines a birational morphism contracting precisely the $(-2)$-curves. {The image is a Du Val cubic surface in $\PP^3$, respectively a degree-four Du Val del Pezzo surface in $\PP^4$; see \cite[\S2.3, Proposition 2.6]{BH}. The degree-four model is cut out by two quadrics; see \cite[Theorem~8.6.2 and Corollary~8.3.7]{DolgachevII}.} Throughout this section, degrees and lines on these images refer to their anticanonical embeddings. We denote the embedded tangent plane at $q\in S^{\sm}$ by $\mathbb T_qS$. The space $\PP(T_aS)$ parametrizes one-dimensional subspaces of $T_aS$, that is, tangent directions at $a$.

These surfaces have only finitely many lines. Indeed, on the minimal resolution, the strict transform $\widetilde\lambda$ of a line has anticanonical degree one. The Hodge index theorem gives $\widetilde\lambda^2\le0$, and adjunction then gives $\widetilde\lambda^2=-1$ and $p_a(\widetilde\lambda)=0$. There are only finitely many such lattice classes, {since the} orthogonal complement of the canonical class in the intersection lattice is negative definite, and a class of negative square contains at most one irreducible curve.

\Needspace{9\baselineskip}
\begin{lemma}\label{lem:centers}
Let $S\subset\PP^4$ be a degree-four Du Val del Pezzo surface in its anticanonical embedding. Fix $a\in S^{\sm}$ and a finite set $Z\subset S\setminus\{a\}$. For $q$ in a nonempty Zariski open subset of $S^{\sm}$, the point $q$ differs from $a$, lies on no line of $S$, and satisfies $a\notin\mathbb T_qS$. The hyperplane $\langle a,\mathbb T_qS\rangle$ cuts an integral curve $C(q,a)$ with an ordinary node at $q$, smooth elsewhere, and disjoint from $Z\cup\Sing S$. Given a direction $\delta\in\PP(T_aS)$, the open subset can also be chosen so that $T_aC(q,a)$ has direction different from $\delta$.
\end{lemma}
\begin{proof}

Let $r:\widetilde S\to S$ be the minimal resolution, let $\beta_a:W\to\widetilde S$ blow up the point above $a$, and set $\mu=r\circ\beta_a$. Write $e_a$ for the exceptional curve of $\beta_a$ and $H=\mu^*(-K_S)$. Crepancy of $r$ and the blow-up formula give
\[
 -K_W=H-e_a,\qquad (-K_W)^2=3.
\]
For an integral curve $C\subset S$ with strict transform $\widehat C\subset W$,
\[
 (-K_W)\cdot\widehat C=\deg C-\operatorname{mult}_a C\ge0.
\]
The inequality follows by intersecting $C$ with a general hyperplane through $a$. The strict transforms of the exceptional curves of $r$ have anticanonical degree zero: $H$ has degree zero on them and they are disjoint from $e_a$. Also $H\cdot e_a=0$ and $e_a^2=-1$, so $(-K_W)\cdot e_a=1$. Hence $W$ is a weak del Pezzo surface of degree three. Its anticanonical morphism
\[
 S\xleftarrow{\ \mu\ }W\xrightarrow{\ \nu\ }S'\subset\PP^3
\]
resolves projection from $a$ and has a Du Val cubic as image. It maps $e_a\simeq\PP(T_aS)$ isomorphically to a line $\ell$. The exceptional curves of $r$, and the strict transforms of any lines of $S$ through $a$, are contracted by $\nu$ to singular points of $S'$.

We use the following two properties of tangent planes to $S'$. First, a general tangent-plane section is an integral cubic with a single ordinary node at its point of tangency. 
 Indeed, the Gauss map is generically finite: otherwise its general fibers, whose closures are linear spaces in characteristic zero \cite[\S1]{MT}, would contradict finiteness of the lines on $S'$. Thus its differential has rank two at a general smooth point $q'$. To interpret this condition, choose local analytic coordinates centered at $q'$, with tangent plane $z=0$, and write the surface as
\[
 z=g_2(\xi_1,\xi_2)+\text{terms of order at least three},
\]
where $g_2$ is homogeneous of degree two. Writing $g=g_2+O(3)$, the two slope coordinates of the Gauss map are $-g_{\xi_1},-g_{\xi_2}$, and the differential of its remaining affine coordinate vanishes at $q'$. Its rank is therefore the rank of the Hessian of $g_2$. Rank two makes $g_2$ nondegenerate. Setting $z=0$ gives a plane curve with an ordinary node at $q'$. Choose $q'$ off the finitely many lines on $S'$. A nonintegral plane cubic has a line component, and any singular point away from its line components would have to be a singular point of an irreducible conic. The tangent-plane section is therefore integral. Since its arithmetic genus is one, its ordinary node at $q'$ is its only singularity.
Second, a general tangent plane avoids any prescribed finite set of points of $\PP^3$. 

If all tangent planes contained a fixed point $p$, the polar derivative in direction $p$ of the cubic equation would vanish on $S'$. Being quadratic, it would vanish identically. The cubic would then be a cone with vertex $p$, of multiplicity three there, contrary to its Du Val singularities. Avoiding $p$ is consequently a nonempty open condition.

\Needspace{4\baselineskip}

{The set $\nu(\mu^{-1}(Z\cup\Sing S))$ is finite, since $Z$ avoids $a$ and $\nu$ contracts every exceptional curve of $r$. Choose a smooth point $q'$ of $S'$, lying on none of its lines, such that its tangent plane $\Pi=\mathbb T_{q'}S'$ cuts an integral cubic whose only singularity is an ordinary node at $q'$. We also require $\Pi$ to avoid $\Sing S'$ and this finite set. If a direction $\delta$ has been prescribed, $\Pi$ must also avoid the corresponding point of $\ell$. Finally, the point $q$ corresponding to $q'$ under projection must belong to $S^{\sm}\setminus\{a\}$ and lie on no line of $S$. Each requirement holds on a nonempty open subset of $S'^{\sm}$. Since $S'$ is irreducible, these requirements can be satisfied simultaneously.}

The cubic $\Gamma=S'\cap\Pi$ avoids the singularities of $S'$, so its inverse image $\widetilde\Gamma\subset W$ avoids the exceptional locus of $\nu$. {Since $\Gamma$ is integral, $\Pi$ does not contain $\ell$. Hence $\Pi\cap\ell$ is a single reduced point, and $\Gamma\cdot\ell=1$.} Since $q'\notin\ell$ and $\Pi$ avoids $\Sing S'$, this intersection is transverse at a smooth point of $\Gamma$ and $S'$. Thus $\widetilde\Gamma$ meets $e_a$ once transversely. Since $\widetilde\Gamma\sim H-e_a$, its image under $\mu$ is an anticanonical curve on $S$, smooth at $a$, with its only singularity an ordinary node at $q$. It avoids $Z\cup\Sing S$, and its tangent direction at $a$ is the point $\Pi\cap\ell$, different from ${{\delta}}$ when prescribed.

Finally, $a\notin\mathbb T_qS$. Otherwise the line $aq$ would meet each defining quadric of $S$ to order at least two at $q$ and also at $a$, so it would lie on both quadrics, contrary to the choice of $q$. Hence $\langle a,\mathbb T_qS\rangle$ is a hyperplane. The hyperplane defining the curve just constructed contains $a$ and is tangent to $S$ at $q$, so it equals this span. The open subset of admissible $q'$ corresponds to a nonempty open subset of $S^{\sm}$, as required.
\end{proof}

\section{Sprays relative to a del Pezzo fibration}

Let $f:V\to B$ be a flat projective morphism from a smooth projective threefold to a smooth projective curve. Assume that every fiber is a degree-four del Pezzo surface with Du Val singularities and that the general fiber is smooth. {Write $V_c=f^{-1}(c)$ for its scheme-theoretic fiber and $V_{B_0}=f^{-1}(B_0)$ for its restriction over an open subset $B_0\subset B$. The relative smooth locus $V^{\sm/B}$ meets each fiber $V_c$ in $V_c^{\sm}$. A local spray $s$ is $f$-relative if $f(s(u,t))=f(u)$.}

\begin{proposition}\label{prop:relative}
Under the hypotheses on $f:V\to B$, at every $a\in V^{\sm/B}$ there are two local $f$-relative rank-one sprays with values in $V$ whose velocities at $a$ span $\ker(df)_a$.
\end{proposition}

{To form the blow-up, we need a section of $f$ through the auxiliary point $q_0$. Smoothness guarantees an \'etale-local section, whereas we need a section over a Zariski neighborhood in $B$.}

\Needspace{5\baselineskip}
\begin{lemma}\label{lem:section}
For every $b\in B$ and $q_0\in V_b^{\sm}$ there is a section $\gamma:B\to V^{\sm/B}$ with $\gamma(b)=q_0$.
\end{lemma}
\begin{proof}
The general fiber is a smooth del Pezzo surface, hence rationally connected. The theorem of Graber--Harris--Starr \cite[Theorem 1.1]{GHS} gives a section $\gamma_*:B\to V$. It lies in $V^{\sm/B}$: differentiating $f\circ \gamma_*=\operatorname{id}_B$ gives a surjective differential of $f$ along it, and both $V$ and $B$ are smooth.

Recall that a smooth variety $Y$ is strongly rationally connected if every point lies on a very free rational curve in $Y$. Here very free means that the curve is given by a morphism $\alpha:\PP^1\to Y$ for which $\alpha^*T_Y$ is ample. 

{Over a smooth projective curve, Theorem 15 of Hassett--Tschinkel \cite{HT} shows that a smooth morphism with strongly rationally connected fibers and a section has a section through any prescribed point of a fiber. The morphism need not be proper.}

The morphism $V^{\sm/B}\to B$ is smooth and has the section $\gamma_*$. Its fiber over $c$ is $V_c^{\sm}$. The surface $V_c$ is a log del Pezzo surface: its anticanonical divisor is ample and its Du Val singularities are Kawamata log terminal. Xu's theorem \cite[Theorem 1.4]{Xu} therefore makes $V_c^{\sm}$ strongly rationally connected. The approximation theorem now gives a section $\gamma:B\to V^{\sm/B}$ with $\gamma(b)=q_0$.
\end{proof}

Blowing up the auxiliary section replaces the degree-four fibers by cubic surfaces with an exceptional line. Planes through that line cut out the conics used for the sprays. 

{For the next lemma, let $g:M\to B'$ be a flat projective morphism from a smooth threefold to a smooth curve whose fibers are Du Val del Pezzo surfaces of fixed degree $d\in\{3,4\}$. Set
\[
 \mathcal A=\omega_{M/B'}^{-1}
 =\omega_M^{-1}\otimes g^*\omega_{B'},
 \qquad \mathcal F=g_*\mathcal A.
\]
Relative projective bundles parametrize one-dimensional quotients; our convention is
\[
 \PP_{B'}(\mathcal F)=\operatorname{Proj}_{B'}(\operatorname{Sym}\mathcal F).
\]}

\begin{lemma}[Relative anticanonical embeddings]\label{lem:relative-embedding}
{The sheaf $\mathcal F$ is locally free of rank $d+1$, and its formation commutes with base change. Evaluation defines a closed immersion
\[
 \iota:M\hookrightarrow\PP_{B'}(\mathcal F),\qquad
 \iota^*\OO_{\PP_{B'}(\mathcal F)}(1)\simeq\mathcal A.
\]
For $d=3$, its image is locally over $B'$ defined by one homogeneous cubic equation with regular coefficients.}
\end{lemma}
\begin{proof}
{Write $\mathcal P=\PP_{B'}(\mathcal F)$.}
Adjunction identifies $\mathcal A|_{{{{M}}}_c}$ with $\omega_{{{{M}}}_c}^{-1}$ for every fiber, since ${{{M}}}_c$ is a Cartier divisor in ${{{M}}}$. This restriction is very ample by the anticanonical description recalled in Section 2. To compute its cohomology, fix a fiber $S$ and let $r:\widetilde S\to S$ be its minimal resolution. Crepancy gives
\[
 K_{\widetilde S}=r^*K_S,\qquad K_{\widetilde S}^2=d,
\]
so $-K_{\widetilde S}$ is nef and big, being the pullback of the ample divisor $-K_S$.

We use the following form of Kawamata--Viehweg vanishing: if $Y$ is smooth and projective over $\CC$ and $D$ is a nef and big Cartier divisor, then $H^i(Y,\OO_Y(K_Y+D))=0$ for $i>0$; see \cite[Theorem 0.1]{FujinoKV}. Applying it on $\widetilde S$ with $D=-K_{\widetilde S}$ and $D=-2K_{\widetilde S}$ gives, respectively,
\[
 H^i(\widetilde S,\OO_{\widetilde S})=0,
 \qquad H^i(\widetilde S,\OO_{\widetilde S}(-K_{\widetilde S}))=0
 \qquad(i>0).
\]
The first vanishing gives $\chi(\OO_{\widetilde S})=1$. Riemann--Roch therefore yields
\[
 h^0(\widetilde S,\OO_{\widetilde S}(-K_{\widetilde S}))
 =1+\tfrac12(-K_{\widetilde S})\cdot(-2K_{\widetilde S})=d+1.
\]
Du Val singularities are rational, so $r_*\OO_{\widetilde S}=\OO_S$ and $R^jr_*\OO_{\widetilde S}=0$ for $j>0$. The projection formula and Leray \cite[Lemmas 20.13.4 and 20.13.6(1)]{StacksLeray} give
\[
 H^i(S,\OO_S(-K_S))\simeq
 H^i(\widetilde S,\OO_{\widetilde S}(-K_{\widetilde S}))\qquad(i\ge0).
\]
Thus $h^0(S,\OO_S(-K_S))=d+1$ and its higher cohomology vanishes.

Since $g$ is projective and $\mathcal A$ is flat over $B'$, cohomology and base change \cite[Lemma 30.22.1]{StacksBaseChange} give the assertions about $\mathcal F$. Evaluation $g^*\mathcal F\to\mathcal A$ is surjective on each fiber, hence surjective by Nakayama's lemma. It defines ${{\iota}}$, whose fiber maps are the complete anticanonical embeddings. The morphism ${{\iota}}$ is proper and quasi-finite, hence finite. Direct image under a finite morphism commutes with base change, so the cokernel of $\mathcal O_{\mathcal P}\to {{\iota}}_*\mathcal O_{{{M}}}$ restricts to zero on every fiber. Nakayama's lemma makes it zero, proving that ${{\iota}}$ is a closed immersion.

For $d=3$, let $p:\mathcal P\to B'$ be the projection and $\mathcal J$ the ideal sheaf of ${{\iota}}({{{M}}})$. In
\[
 0\longrightarrow\mathcal J\longrightarrow\mathcal O_{\mathcal P}
 \longrightarrow {{\iota}}_*\mathcal O_{{{M}}}\longrightarrow0,
\]
both the middle and right-hand terms are flat over $B'$. Thus $\mathcal J$ is flat over $B'$, and restriction to a fiber remains exact. Since ${{\iota}}({{{M}}}_c)$ is a cubic hypersurface, $\mathcal J|_{\mathcal P_c}\simeq\mathcal O_{\mathbf P^3}(-3)$. Cohomology and base change applied to $\mathcal J(3)$ show that ${{\mathcal L_{\mathrm{cub}}}}=p_*\mathcal J(3)$ is a line bundle. Evaluation $p^*{{\mathcal L_{\mathrm{cub}}}}\to\mathcal J(3)$ is an isomorphism on every fiber. Its cokernel vanishes by Nakayama; its kernel also vanishes by Nakayama, since flatness of $\mathcal J(3)$ preserves exactness on restriction to a fiber. Consequently,
\[
 \mathcal J\simeq p^*{{\mathcal L_{\mathrm{cub}}}}\otimes\mathcal O_{\mathcal P}(-3).
\]
After trivializing $\mathcal F$ and ${{\mathcal L_{\mathrm{cub}}}}$ over an open subset of $B'$, the inclusion $\mathcal J\hookrightarrow\mathcal O_{\mathcal P}$ is given by one homogeneous cubic polynomial with regular coefficients there.
\end{proof}

{Let $\gamma:B\to V^{\sm/B}$ be a section, and let $b\in B$ be such that $\gamma(b)$ lies on no line of $V_b$. For each open neighborhood $B_0\subset B$ of $b$, denote the blow-up by
\[
 \beta:M=\Bl_{\gamma(B_0)}V_{B_0}\longrightarrow V_{B_0},
\]
set $g=f|_{V_{B_0}}\circ\beta$, and let $E\subset M$ be its exceptional divisor.}

\begin{lemma}[{Blow-up along the auxiliary section}]\label{lem:section-blowup}
{There is a Zariski open neighborhood $B_0$ of $b$ for which $M$ is smooth and $g$ is flat and projective. For every $c\in B_0$, its scheme-theoretic fiber is
\[
 M_c\simeq\Bl_{\gamma(c)}V_c,
\]
a degree-three Du Val del Pezzo surface. The exceptional divisor $E\to B_0$ is a $\PP^1$-bundle. Each fiber $E_c$ is a line in the smooth locus of the anticanonically embedded cubic surface $M_c$. Moreover,
\[
 \omega_M^{-1}\otimes g^*\omega_{B_0}
 \simeq\beta^*(\omega_{V/B}^{-1}|_{V_{B_0}})\otimes\OO_M(-E).
\]}
\end{lemma}
\begin{proof}
By Lemma~\ref{lem:relative-embedding}, $V$ is anticanonically embedded in a $\PP^4$-bundle.

{The locus of points lying on lines of their fibers is closed: it is the proper image of the universal family over the relative Hilbert scheme of lines. Since $\gamma(b)$ is outside this locus, $\gamma(c)$ lies on no line of $V_c$ for $c$ in some open neighborhood $B_0$ of $b$.}

{The center $Z_\gamma=\gamma(B_0)$ is a smooth curve in the smooth threefold $V_{B_0}$, so its blow-up $M$ is smooth.} Along $Z_\gamma$, the morphism $f$ is smooth of relative dimension two. \'Etale locally, two relative coordinates $x,y$ generate the ideal of $Z_\gamma$ and identify the blow-up with the \'etale pullback of the blow-up of $\AA^2_{B_0}$ along its zero section. This model is flat over $B_0$ and commutes with taking fibers. Away from $Z_\gamma$, the blow-up agrees with $V_{B_0}$. Thus $g$ is flat and $M_c\simeq\Bl_{\gamma(c)}V_c$ scheme-theoretically. 
If $\mathcal I_{Z_\gamma}$ is the ideal of $Z_\gamma$, then $E\simeq\mathbf P_{Z_\gamma}(\mathcal I_{Z_\gamma}/\mathcal I_{Z_\gamma}^2)$; the conormal bundle has rank two, so $E\to B_0$ is a $\PP^1$-bundle.

Fix $c\in B_0$, set $q_c=\gamma(c)$, and write $e_c=E|_{M_c}$. Let $H$ be the pullback of $-K_{V_c}$ to $M_c$. Since $q_c$ is a smooth point,
\[
 -K_{M_c}=H-e_c,\qquad (H-e_c)^2=3,\qquad (H-e_c)\cdot e_c=1.
\]
On the strict transform of an integral curve $D\subset V_c$, the divisor $H-e_c$ has degree $\deg D-\operatorname{mult}_{q_c}D$. This is nonnegative, as is seen by intersecting $D$ with a general hyperplane through $q_c$. If it were zero, projection from $q_c$ would contract $D$. A fiber of that projection lies on a line through $q_c$, so $D$ itself would be such a line, contrary to the choice of the section. Together with the positive degree on $e_c$ and the positive self-intersection, the Nakai--Moishezon criterion shows that $H-e_c$ is ample.

The singularities of $M_c$ are the Du Val singularities of $V_c$ away from $q_c$; the blow-up is smooth near $e_c$. Thus $M_c$ is a degree-three Du Val del Pezzo surface. The canonical-bundle formula for a blow-up with smooth center of codimension two gives the asserted formula on $M$. Its restriction to $M_c$ is $\mathcal O_{M_c}(H-e_c)$. By Lemma~\ref{lem:relative-embedding}, $M$ is embedded as a family of cubic surfaces. Since $(H-e_c)\cdot e_c=1$, each $e_c$ is embedded as a line in the smooth locus of its cubic.
\end{proof}

\begin{lemma}[Conic pencils on cubic surfaces]\label{lem:cubic-conics}
Let $g:M\to B'$ be a flat projective family of degree-three Du Val del Pezzo surfaces, where $M$ is a smooth threefold and $B'$ is a smooth curve. Suppose that $E\subset M$ is a $\PP^1$-bundle over $B'$ whose fibers are lines in the anticanonical embeddings and lie in the smooth loci of the surface fibers. The planes containing these lines form a $\PP^1$-bundle $\tau:T\to B'$. Their residual conics define a proper flat morphism $h:M\to T$, with $\tau\circ h=g$, whose scheme-theoretic fibers are those conics.
\end{lemma}
\begin{proof}
By Lemma~\ref{lem:relative-embedding}, $M$ is embedded in ${\mathbf P_{B'}}(\mathcal F)$, where $\mathcal F=g_*\omega_{M/B'}^{-1}$ has rank four. We write $\mathcal O_M(1)$ for the hyperplane bundle of this embedding; thus $\mathcal O_M(1)\simeq\omega_{M/B'}^{-1}$. Restriction of linear forms to $E_c$ is surjective. Since $\mathcal O_M(1)|_{E_c}\simeq\mathcal O_{\PP^1}(1)$, base change gives a surjection
\[
 \mathcal F\longrightarrow (g|_E)_*\mathcal O_E(1)
\]
onto a rank-two bundle. Let $\mathcal K$ be its rank-two kernel. A plane containing $E_c$ is specified by a nonzero linear form in $\mathcal K_c$, up to scalar. With our quotient convention for projectivization, these pencils therefore form
\[
 T={\mathbf P_{B'}}(\mathcal K^\vee)\longrightarrow B'.
\]

A plane through $E_c$ cuts out that line and a residual conic. We check that the residual pencil defines a morphism even along $E$. Locally on $B'$, choose relative projective coordinates so that
\[
 M=\{x\Psi_1(x,y,z,w)+y\Psi_2(x,y,z,w)=0\},\qquad E=\{x=y=0\},
\]
where $\Psi_1,\Psi_2$ are quadratic forms with coefficients on the base. The existence of a relative cubic equation follows from Lemma~\ref{lem:relative-embedding}; coordinates adapted to $E$ come from locally splitting the preceding surjection. Off $E$, the plane containing a point and $E_c$ is recorded by $[x:y]$. The cubic equation gives
\[
 [x:y]=[-\Psi_2:\Psi_1]
\]
where both expressions are defined. On $E_c$, the partial derivatives of the cubic equation are $(\Psi_1,\Psi_2,0,0)$. Smoothness of $M_c$ along $E_c$ says that $\Psi_1,\Psi_2$ have no common zero there. The second expression extends the map across $E$; the local extensions agree off $E$ and hence everywhere, giving $h:M\to T$.

For $[\lambda:\mu]\in T_c$, the corresponding plane is $\mu x-\lambda y=0$. Substituting $x=\lambda r$, $y=\mu r$ in the cubic equation gives
\[
 r\bigl(\lambda \Psi_1(\lambda r,\mu r,z,w)
       +\mu \Psi_2(\lambda r,\mu r,z,w)\bigr)=0.
\]
The factor $r=0$ defines $E_c$. The quadratic factor defines the residual conic and is not identically zero, since an integral cubic surface contains no plane. It cannot contain $E_c$: that would make $\lambda \Psi_1+\mu \Psi_2$ vanish identically on $E_c$, forcing the two quadratic forms $\Psi_1|_{E_c},\Psi_2|_{E_c}$ to have a common zero on $\PP^1$, contrary to smoothness along $E_c$. Off $E_c$ it is the scheme-theoretic fiber of $h_c$. Near a point of $E_c$ where $\Psi_1\ne0$, normalize $\mu=1$. The fiber equation is $\Psi_2+\lambda \Psi_1=0$; together with the cubic equation it implies $x=\lambda y$. Conversely, in that plane the residual equation is $\lambda \Psi_1+\Psi_2=0$. This proves the same scheme-theoretic identification along $E_c$; where $\Psi_2\ne0$ the argument is symmetric. Thus every fiber of $h$ is a plane conic of pure dimension one, possibly reducible or nonreduced.

The morphism $h$ is proper because $M$ is projective over $B'$ and $T$ is separated over $B'$. The total spaces $M$ and $T$ are smooth, of dimensions three and two. The criterion called \emph{miracle flatness} states that a local homomorphism $R_{\mathrm{base}}\to R_{\mathrm{tot}}$ of Noetherian local rings is flat if $R_{\mathrm{base}}$ is regular, $R_{\mathrm{tot}}$ is Cohen--Macaulay, and
\[
 \dim R_{\mathrm{tot}}=\dim R_{\mathrm{base}}+\dim(R_{\mathrm{tot}}/\mathfrak m_{R_{\mathrm{base}}} R_{\mathrm{tot}});
\]
see \cite[Lemma 10.128.1, Tag 00R4]{Stacks}. Here $R_{\mathrm{tot}}/\mathfrak m_{R_{\mathrm{base}}} R_{\mathrm{tot}}$ is the local ring of the scheme-theoretic fiber.

{At a closed point $m\in M$, with $\theta=h(m)$, the relevant local rings are}
\[
 R_{\mathrm{base}}=\mathcal O_{T,{\theta}},\qquad R_{\mathrm{tot}}=\mathcal O_{M,{m}}.
\]
The point ${\theta}$ is closed, since we work with varieties over $\CC$. Smoothness of $T$ makes $R_{\mathrm{base}}$ regular of dimension two; smoothness of $M$ makes $R_{\mathrm{tot}}$ regular, hence Cohen--Macaulay, of dimension three. The fiber $M_{\theta}$ is a plane conic of pure dimension one, so its local ring at $m$ has dimension one. Thus the required equality is $3=2+1$, and $h$ is flat at every closed point of $M$. The flat locus is open for this morphism of finite presentation between Noetherian schemes. Its complement, if nonempty, would contain a closed point of the complex variety $M$. Hence $h$ is flat everywhere.
\end{proof}

\Needspace{22\baselineskip}

\begin{proposition}[Relative conic pencil]\label{lem:conic}
{Let $\gamma:B\to V^{\sm/B}$ be a section, and let $b\in B$ be such that $\gamma(b)$ lies on no line of $V_b$. There is a Zariski open neighborhood $B_0$ of $b$ such that the blow-up
\[
 \beta:M=\Bl_{\gamma(B_0)}V_{B_0}\longrightarrow V_{B_0}
\]
admits a flat morphism $h:M\to T$ whose fibers are plane conics, where $\tau:T\to B_0$ is a $\PP^1$-bundle. With $g=f|_{V_{B_0}}\circ\beta$, the diagram}
\[
\begin{tikzcd}[column sep=large,row sep=large]
 M \arrow[r,"h"] \arrow[dr,"g"'] & T \arrow[d,"\tau"] \\
 & B_0
\end{tikzcd}
\]
{commutes. For each $c\in B_0$, the restriction
\[
 h_c:M_c\longrightarrow T_c\simeq\PP^1
\]
is the conic pencil obtained from the anticanonical curves of $V_c$ singular at $\gamma(c)$ by blowing up that point.}
\end{proposition}
\begin{proof}
Lemma~\ref{lem:section-blowup} gives $B_0$, the blow-up $\beta:M\to V_{B_0}$, and its exceptional divisor $E$. The fibers of $g$ are cubic del Pezzo surfaces, and $E_c$ is a line in the smooth locus of $M_c$. Lemma~\ref{lem:cubic-conics} gives the bundle $T\to B_0$ and the proper flat conic morphism $h:M\to T$.

{To identify the pencil on the original surface, for $c\in B_0$ set $q_c=\gamma(c)$, $e_c=E_c$, and $H=\beta_c^*(-K_{V_c})$.} A plane through $e_c$ cuts out $e_c$ and a residual conic of class
\[
 (-K_{M_c})-e_c=(H-e_c)-e_c=H-2e_c.
\]
The planes through $e_c$ give the complete system $|-K_{M_c}-e_c|$, since the cubic embedding is anticanonical. Under blow-up, sections of $\mathcal O_{M_c}(H-2e_c)$ correspond to anticanonical sections on $V_c$ vanishing to order at least two at $q_c$. These are precisely the anticanonical curves singular there. No residual conic contains $e_c$, as shown in Lemma~\ref{lem:cubic-conics}, so its image has multiplicity two at $q_c$ and the conic is its strict transform. For the member chosen in Lemma~\ref{lem:centers}, the blow-up resolves its ordinary node.
\end{proof}

\begin{proof}[Proof of Proposition~\ref{prop:relative}]
{Set $b=f(a)$ and $S=V_b$. Lemma~\ref{lem:centers} supplies an auxiliary point $q_0$, with $Z=\varnothing$, and 
{Lemma~\ref{lem:section} gives a section $\gamma$ with $\gamma(b)=q_0$.} Proposition~\ref{lem:conic} gives the relative conic pencil $h:M\to T$.} The strict transform $\widetilde C$ of $C(q_0,a)$ in $M_b$ is a smooth complete conic: blowing up resolves its ordinary node at $q_0$, and the curve has no other singularities and avoids $\Sing S$.

Let ${{\theta_0}}\in T_b$ be its conic parameter. Since $h$ is proper and flat, the parameters of nonsmooth fibers form a closed subset of $T$. There is therefore an open neighborhood $T^\circ$ of ${{\theta_0}}$ over which
\[
 h^\circ:M^\circ=h^{-1}(T^\circ)\longrightarrow T^\circ
\]
is a smooth proper conic family. Choose an open neighborhood $U$ of $a$ in $V_{B_0}$, disjoint from $\gamma(B_0)$, whose lift $\eta:U\to M$ lies in $M^\circ$. 
{Pullback along $h\circ\eta$ gives the Cartesian square}
\[
\begin{CD}
 \widehat U=M^\circ\times_{T^\circ}U @>{\pr_M}>> M^\circ\\
 @V{{{p_U}}}VV @VV{h^\circ}V\\
 U @>>{h\circ\eta}> T^\circ.
\end{CD}
\]

{The fiber of $p_U$ over $u$ is the smooth conic through $\eta(u)$. Thus the pullback has the section}
\[
 \sigma:U\longrightarrow \widehat U,\qquad \sigma(u)=(\eta(u),u).
\]

The image $\sigma(U)$ is an effective Cartier divisor on $\widehat U$, and $\OO_{\widehat U}(\sigma(U))$ has degree one on every fiber. Its direct image under ${{p_U}}$ is locally free of rank two, has formation compatible with base change, and its evaluation map realizes $\widehat U$ as the associated projective bundle. These assertions follow from $h^0(\PP^1,\OO(1))=2$, $h^1(\PP^1,\OO(1))=0$, and the fiberwise embedding by $\OO(1)$. {After shrinking $U$ around $a$, there is a trivialization}
\[
 \widehat U\simeq U\times\PP^1
\]
carrying $\sigma$ to the constant section $0$. Such a trivialization exists by locally splitting the rank-one quotient associated with the section. Its affine chart gives an open immersion $j:U\times\AA^1\hookrightarrow \widehat U$ with $j(u,0)=\sigma(u)$.

{The resulting local spray is}
\begin{equation}\label{eq:spray}
 s_{q_0}=\beta\circ\pr_M\circ j:U\times\AA^1\longrightarrow V.
\end{equation}

{This morphism satisfies $s_{q_0}(u,0)=\beta(\eta(u))=u$.} The projections from $\widehat U$ commute with the maps to $B_0$, so $f(s_{q_0}(u,t))=f(u)$. On the fiber over $a$, the map $j$ parametrizes an affine open subset of $\widetilde C$. Since $\beta$ is an isomorphism near $\eta(a)$, the velocity of $s_{q_0}$ at $a$ is nonzero and spans $T_aC(q_0,a)$.

{Lemma~\ref{lem:centers}, with $\delta$ the direction of $T_aC(q_0,a)$, supplies a second auxiliary point. The same construction gives a second nonzero velocity in $T_aS$ with a different direction.} These two velocities span $T_aS=\ker(df)_a$. Intersecting their base neighborhoods gives two local sprays at the same point, both with values in $V$.
\end{proof}

\begin{corollary}\label{cor:criterion}
Let a smooth projective threefold $V$ admit finitely many fibrations $f_i:V\to B_i$ satisfying the hypotheses of Proposition~\ref{prop:relative}. Suppose that at every $a\in V$ there are indices $i,j$ such that $(f_i,f_j):V\to B_i\times B_j$ has differential of rank two at $a$. Then $V$ is algebraically elliptic.
\end{corollary}
\begin{proof}
{Fix $a\in V$ and choose $i,j$ as in the hypothesis. Since $d(f_i,f_j)_a$ has rank two, both $(df_i)_a$ and $(df_j)_a$ are surjective. Thus $f_i$ and $f_j$ are smooth at $a$, and the tangent planes to their fibers are
\[
 K_i:=\ker(df_i)_a,\qquad K_j:=\ker(df_j)_a.
\]
The rank-two hypothesis also implies that these planes are distinct. Being two-dimensional subspaces of the three-dimensional space $T_aV$, they satisfy $K_i+K_j=T_aV$.}

Proposition~\ref{prop:relative} supplies two local $f_i$-relative rank-one sprays whose velocities at $a$ span $K_i$, and two local $f_j$-relative rank-one sprays whose velocities span $K_j$. Restricting their base domains to a common Zariski neighborhood of $a$, we obtain four local sprays with values in $V$ whose velocities span $T_aV$. Since this holds at every $a$, the localization principle and the equivalence of algebraic subellipticity and ellipticity recalled in the introduction yield the conclusion.
\end{proof}

\section{The double covers}

{With $Y=(\PP^1)^3$ and $L=\OO_Y(1,1,1)$, let $\rho\in H^0(Y,L^{\otimes2})$ define the smooth branch divisor $R$.} Multiplication by $\rho$ gives a map $L^{-2}\to\OO_Y$ and hence an algebra structure on $\OO_Y\oplus L^{-1}$. The relative affine spectrum provides the associated double cover
\[
 \pi:X=\operatorname{Spec}_Y(\OO_Y\oplus L^{-1})\longrightarrow Y.
\]
On an affine open set $U$ trivializing $L$, with $\rho$ represented by $F\in \OO(U)$, the cover is
\[
 \pi^{-1}(U)=\operatorname{Spec}\bigl(\OO(U)[w]/(w^2-F)\bigr).
\]
Thus $\pi$ is finite flat of degree two, with branch divisor $R$. The variety $X$ is smooth because $R$ is smooth, and projective because $\pi$ is finite and $Y$ is projective. It is integral: along any component of the nonempty reduced divisor $R$, 
the defining rational square class has odd valuation. 

\begin{lemma}\label{lem:fibers}
Every fiber of each $f_i$ is an integral normal del Pezzo surface of degree four with at worst Du Val singularities. The general fiber is smooth.
\end{lemma}
\begin{proof}
Fix a coordinate slice $Q\simeq\PP^1\times\PP^1$ in $Y$. The corresponding fiber $S$ is the double cover of $Q$ branched along the scheme-theoretic intersection $D=R\cap Q$. {We denote this cover by $\pi_S:S\to Q$.} Smoothness of $R$ allows $D$ to be singular. 

Choose a coordinate ${{\zeta}}$ on an affine chart of the omitted factor, with $Q=\{{{\zeta}}=0\}$, and trivialize the corresponding factor of $\OO_Y(2,2,2)$. The branch equation has the form
\[
 F({{\zeta}};u,v)=F_0(u,v)+{{\zeta}} F_1(u,v)+{{\zeta}}^2F_2(u,v),
 \qquad F_j\in H^0(Q,\OO_Q(2,2)).
\]
Along $Q$, the derivative in the transverse direction is $F_1$, whereas the two derivatives tangent to $Q$ are those of $F_0$.

Suppose first that $F_0$ vanishes identically, so $Q\subset R$. Both tangent derivatives then vanish everywhere on $Q$. The section $F_1$ must have a zero: a nowhere-vanishing section would trivialize the nontrivial line bundle $\OO_Q(2,2)$. At such a zero the transverse derivative vanishes as well, making $R$ singular. This contradiction shows that $F_0\ne0$. Its zero divisor $D$ is nonempty because its line bundle is nontrivial.

Suppose that $D$ has a multiple irreducible component $C$. Both tangent derivatives of $F_0$ vanish along $C$. The pullback of $F_1|_C$ to the normalization of $C$ is a section of a line bundle of positive degree, since $\OO_Q(2,2)$ is ample. {This section has a zero.}
 At the image of that zero, $F$ and all three of its derivatives vanish. This again contradicts smoothness of $R$. Thus $D$ is reduced, and its singular points are isolated.

Near a point of $Q$, trivialize the double cover and write it as
\[
 w^2={{G}}(u,v),
\]
where ${{G}}=0$ is a local equation of $D$. The Jacobian criterion says that a singular point of the surface must satisfy
\[
 w=0,\qquad {{G}}={{G}}_u={{G}}_v=0.
\]
The surface is therefore smooth over $Q\setminus\Sing D$. At a singular point of $D$, choose affine product coordinates $u,v$ centered there. Before making further analytic coordinate changes, ${{G}}$ is a polynomial of degree at most two in each variable. Its multiplicity is consequently at most four.

If the multiplicity is two and the quadratic part has rank two, the analytic Morse lemma puts ${{G}}$ in the form $uv$. The surface equation becomes
\[
 w^2=uv,
\]
which is an $A_1$ singularity. If the quadratic part has rank one, the analytic splitting lemma gives ${{G}}=u^2+\varphi(v)$, where $\varphi$ has order at least three. The series $\varphi$ is nonzero, since otherwise $D$ would not be reduced. Write $\varphi(v)=v^m a(v)$, with $m\ge3$ and $a(0)\ne0$. Absorbing an analytic $m$th root of $a(v)$ into $v$ gives ${{G}}=u^2+v^m$. Hence the surface equation is
\[
 (w-u)(w+u)=v^m.
\]
The variables $w-u,w+u,v$ are analytic coordinates, so this is the normal form of an $A_{m-1}$ singularity.

If the multiplicity is three, the bidegree restriction leaves only the monomials $u^2v$, $uv^2$, and $u^2v^2$. Thus
\[
 {{G}}=uv(au+bv+cuv).
\]
Reducedness forces $a\ne0$ and $b\ne0$: if either vanished, $u$ or $v$ would be a repeated factor. The three branches have distinct tangent lines, so this is an ordinary triple point. Rescaling $u,v,w$ puts the surface equation in the form
\[
 w^2=uv(u+v+\kappa uv)
\]
for a constant $\kappa$. {The analytic change of coordinates
\[
 v'=v(1+\kappa u),\qquad w'=w\sqrt{1+\kappa u}
\]
is defined near the origin, since $1+\kappa u$ is a unit. It gives
\[
 (w')^2=uv'(u+v'),
\]
the $D_4$ normal form.} Finally, multiplicity four would force ${{G}}=cu^2v^2$, contrary to reducedness. This exhausts the possibilities and proves that the surface has only Du Val singularities.

{The surface $S$ is integral.} Indeed, at each component of the nonempty reduced branch divisor its defining square class has odd valuation. It is therefore not a square in $\CC(Q)$, so the generic double cover is a field extension of degree two. The locally free algebra defining $S$ injects into that field, proving integrality. Locally $S$ is a hypersurface, so it is Cohen--Macaulay and satisfies Serre's condition $S_2$. Its singular set is finite by the calculation above, so it is regular in codimension one. Serre's criterion gives normality.

{The canonical bundle formula for the double cover gives
\[
 \omega_S\simeq\pi_S^*\bigl(\omega_Q\otimes\OO_Q(1,1)\bigr),\qquad
 \omega_S^{-1}\simeq\pi_S^*\OO_Q(1,1).
\]}
{Thus $-K_S$ is ample because $\pi_S$ is finite.} To compute its self-intersection, let $D_1$ and $D_2$ be the divisor classes of fibers of the two projections $Q\to\PP^1$. 

\[
 D_1^2=D_2^2=0,\qquad D_1\cdot D_2=1,\qquad
 c_1(\OO_Q(1,1))=D_1+D_2.
\]
Writing ${{H_Q}}=D_1+D_2$, we obtain ${{H_Q}}^2=2$. Since $\pi_S$ has degree two, the projection formula for intersection numbers gives
\[
 K_S^2=(-K_S)^2=(\pi_S^*{{H_Q}})^2
 =\deg(\pi_S)\,{{H_Q}}^2
 =2(D_1+D_2)^2=2\cdot2=4.
\]

Thus $S$ is a degree-four del Pezzo surface. The morphism $f_i$ is flat, as the composition of a finite flat cover and a projection. Bertini's theorem applied to its base point free pencil on the smooth threefold $X$ shows that the general fiber is smooth.
\end{proof}

\begin{proof}[Proof of Theorem~\ref{thm:main}]
Lemma~\ref{lem:fibers} verifies the fiber hypotheses for all three $f_i$. We check the rank condition in Corollary~\ref{cor:criterion}.

Off the ramification divisor of $\pi$, the differential $d\pi$ has rank three, so the three coordinate differentials are independent. At a ramification point $x\in X$, use local coordinates $u_1,u_2,u_3$ on $Y$ and a local trivialization of $L$. The equation is $w^2=F(u_1,u_2,u_3)$, with $w=0$ at the point. Set $y=\pi(x)$. In the chosen local trivialization, write an ambient tangent vector at $x$ as $(a,v)\in\CC\oplus T_yY$, where $a$ is its component in the fiber coordinate $w$. Differentiating the defining equation gives $2w(x)a-dF_y(v)=0$. Since $w(x)=0$, this identifies
\[
 T_xX=\{(a,v)\in\CC\oplus T_yY:dF_y(v)=0\},
 \qquad d\pi_x(a,v)=v.
\]
Smoothness of $R$ gives $dF_y\ne0$, so the image of $d\pi_x$ is the two-dimensional plane $\ker(dF_y)$. The restrictions of $du_1,du_2,du_3$ span its dual, so some pair is independent. Their pullbacks to $T_xX$ are the differentials of the corresponding $f_i$ and remain independent. Corollary~\ref{cor:criterion} therefore applies and proves the theorem.
\end{proof}

\Needspace{16\baselineskip}
\begin{proposition}\label{prop:unirational}
The threefold $X$ admits a dominant rational map $\PP^3\dashrightarrow X$ of degree two.
\end{proposition}
\begin{proof}
Set $Q=\PP^1\times\PP^1$ and ${{\varpi}}=(f_1,f_2):X\to Q$. Its generic fiber is a smooth conic: it is the double cover of the third $\PP^1$ branched at two points. Choose a general $c\in\PP^1$ and set $S=f_3^{-1}(c)$. This is a smooth degree-four del Pezzo surface, hence rational, and the restriction ${{\varpi}}|_S:S\to Q$ is finite of degree two.

The base change of the conic fibration ${{\varpi}}:X\to Q$ along ${{\varpi}}|_S:S\to Q$ is the projection $X\times_Q S\to S$. {Since the generic conic of $\varpi$ is geometrically integral, the base change has a unique component $W$ dominating $S$. The inclusion $S\hookrightarrow X$ gives a section $s\mapsto(s,s)$ of $W\to S$. Thus its generic conic has a $\CC(S)$-rational point and is isomorphic to $\PP^1_{\CC(S)}$. Consequently, $W$ is birational to $S\times\PP^1$ and is rational. Geometric integrality also implies that $\CC(X)\otimes_{\CC(Q)}\CC(S)$ is a field, of degree two over $\CC(X)$; hence the projection $W\to X$ has degree two.} A birational parametrization by $\PP^3$ gives the assertion.
\end{proof}


\begin{thebibliography}{99}
\footnotesize
\setlength{\itemsep}{0pt}
\setlength{\parsep}{0pt}
\bibitem{Ban} J. Banecki,
 \emph{Retract rational varieties are uniformly retract rational}, 
arXiv:2411.17892 (2024), revised 2026.
\bibitem{BH} R.~Blache and E.~Hallouin, \emph{Construction of good codes from weak Del Pezzo surfaces}, arXiv:2301.11825 (2023), \S2.3. \url{https://arxiv.org/abs/2301.11825}.
\bibitem{DolgachevII} {I.~V.~Dolgachev, \emph{Classical Algebraic Geometry: a modern view. II}, draft, August 10, 2025, Theorem~8.6.2 and Corollary~8.3.7. \url{https://sites.lsa.umich.edu/idolga/wp-content/uploads/sites/1467/2025/10/CAG-2.pdf}.}
\bibitem{DKZ} A.~Dvorsky, S.~Kaliman,  and M.~Zaidenberg,
\emph{Algebraic ellipticity of quartic double spaces},
  arXiv:2610.07260 (2026).
\bibitem{FujinoKV} O.~Fujino, \emph{Kawamata--Viehweg vanishing theorem}, notes, July 22, 2009, Theorem~0.1. \url{https://www.math.kyoto-u.ac.jp/~fujino/kawamata--viehweg2.pdf}.
\bibitem{GHS} T.~Graber, J.~Harris, and J.~Starr, \emph{Families of rationally connected varieties}, J. Amer. Math. Soc. \textbf{16} (2003), no.~1, 57--67. \url{https://arxiv.org/abs/math/0109220}.
\bibitem{HT} B.~Hassett and Y.~Tschinkel, \emph{Approximation at places of bad reduction for rationally connected varieties}, Pure Appl. Math. Q. \textbf{4} (2008), no.~3, 743--766. \url{https://arxiv.org/abs/math/0605274}.
\bibitem{KZ} S.~Kaliman and M.~Zaidenberg, \emph{Gromov ellipticity and subellipticity}, Forum Math. \textbf{36} (2024), no.~2, 373--376. \url{https://arxiv.org/abs/2301.03058}. {The theorem and lemma numbers cited here refer to arXiv:2301.03058v1.}
\bibitem{KZCub} {S.~Kaliman and M.~Zaidenberg, \emph{Algebraic Gromov's ellipticity of cubic hypersurfaces}, Proc. Steklov Inst. Math. \textbf{329} (2025), 79--87. \url{https://doi.org/10.1134/S008154382560067X}.}
\bibitem{MT} E.~Mezzetti and O.~Tommasi, \emph{Some remarks on varieties with degenerate Gauss image}, Pacific J. Math. \textbf{213} (2004), no.~1, 79--88. \url{https://msp.org/pjm/2004/213-1/pjm-v213-n1-p05-s.pdf}.
\bibitem{MM} S.~Mori and S.~Mukai, \emph{Classification of Fano 3-folds with $B_2\ge2$}, Manuscripta Math. \textbf{36} (1981/82), no.~2, 147--162. \url{https://doi.org/10.1007/BF01170131}.
\bibitem{Prokhorov} Yu.~Prokhorov, \emph{Rationality of Fano threefolds with terminal Gorenstein singularities, I}, arXiv:1907.05678 (2019), Theorem~6.1. \url{https://arxiv.org/abs/1907.05678}.
\bibitem{Stacks} The Stacks Project Authors, \emph{The Stacks Project}, Lemma 10.128.1 (miracle flatness), Tag 00R4. \url{https://stacks.math.columbia.edu/tag/00R4}.
\bibitem{StacksLeray} The Stacks Project Authors, \emph{The Stacks Project}, Cohomology of Sheaves, Lemmas 20.13.4 and 20.13.6(1), Tags 01F2 and 01F4. \url{https://stacks.math.columbia.edu/tag/01F2}; \url{https://stacks.math.columbia.edu/tag/01F4}.
\bibitem{StacksBaseChange} The Stacks Project Authors, \emph{The Stacks Project}, Cohomology of Schemes, Lemma 30.22.1 (cohomology and base change), Tag 07VK. \url{https://stacks.math.columbia.edu/tag/07VK}.
\bibitem{Xu} C.~Xu, \emph{Strong rational connectedness of surfaces}, J. Reine Angew. Math. \textbf{665} (2012), 189--205. \url{https://arxiv.org/abs/0810.2597}.
\end{thebibliography}
\end{document}